%
\input gtmacros
\input amsnames
\input gtmonout
\volumenumber{1}
\volumeyear{1998}
\volumename{The Epstein birthday schrift}
\pagenumbers{413}{418}
\received{15 November 1997}
\published{27 October 1998}
\papernumber{20}

\def\R{{\Bbb R}}              
\def\Z{\hbox{${\Bbb Z}$}}

\def\setminus{-}

\def\deg#1{\hbox{Degree}(#1)}

\reflist

\refkey\BKS
{\bf R\,G~Burns}, {\bf A~Karrass}, {\bf D~Solitar}
{\it A note on groups with separable finitely generated
subgroups},
Bull. Aust. Math. Soc.
36
(1987)
153--160

\refkey\N
{\bf G\,A~Niblo},
{\it Finding splittings of groups and 3--manifolds},
Bull. London Math. Soc.
27
(1995)
567--574

\refkey\Long
{\bf D\,D~Long},
{\it Engulfing and subgroup separability for hyperbolic 3--manifolds},
Trans. Amer. Math. Soc.
308
(1988)
849--859

\refkey\LN
{\bf D\,D Long}, {\bf G\,A Niblo},
{\it Subgroup separability and 3--manifold groups},
Math. Z.
207
(1991)
209--215

\refkey\NW
{\bf G\,A~Niblo}, {\bf D\,T~Wise},
{\it Subgroup separability, knot groups, and graph manifolds},
Proc. Amer. Math. Soc. (to appear)

\refkey\RW
{\bf H Rubinstein}, {\bf S Wang},
{\it  $\pi_1$--injective surfaces in graph manifolds},
Comm. Math. Helv. (to appear)

\refkey\Sc
{\bf G\,P Scott},
{\it Subgroups of surface groups are almost geometric},
J. London Math. Soc.
17
(1978)
555--565

\endreflist

\title{The engulfing property for 3--manifolds}      
\asciititle{The engulfing property for 3-manifolds}      

\authors{Graham A Niblo\\Daniel T Wise}     
\asciiauthors{Graham A Niblo and Daniel T Wise}

\address{Faculty of Mathematical Studies, University of Southampton\\
Highfield, Southampton, SO17 1BJ, UK\\{\rm and}\\
Department of Mathematics, Cornell University\\Ithaca, NY 14853, USA}
\asciiaddress{Faculty of Mathematical Studies, University of Southampton\\
Highfield, Southampton, SO17 1BJ, UK\\
Department of Mathematics, Cornell University\\Ithaca, NY 14853, USA}

\email{gan@maths.soton.ac.uk, daniwise@math.cornell.edu}

\abstract
We  show that there are Haken 3--manifolds whose
fundamental groups do not satisfy the engulfing property. In
particular one can construct a $\pi_1$--injective immersion of a surface into
a
graph manifold which does not factor through any proper finite cover of the
3--manifold.
\endabstract

\asciiabstract{%
We  show that there are Haken 3-manifolds whose
fundamental groups do not satisfy the engulfing property. In
particular one can construct a pi_1-injective immersion of a surface into
a
graph manifold which does not factor through any proper finite cover of the
3-manifold.}

\primaryclass{20E26}\secondaryclass{20F34, 57M05, 57M25}

\keywords{Double coset decompositions, subgroup separability,
3--man\-ifolds, engulfing property}
\asciikeywords{Double coset decompositions, subgroup separability,
3-manifolds, engulfing property}

\maketitle

\section{Introduction}
\rk{Definition} A subgroup $H$ of a group $G$ is said to be
{\bf separable} if it is an intersection of finite index subgroups of $G$.
It is said to be {\bf engulfed} if it is contained in a
proper subgroup of finite index in $G$.

Subgroup separability was first explored as a tool in low dimensional
topology by Scott in [\Sc]. He showed that if $f\co
\Sigma\longrightarrow M$ is a $\pi_1$--injective immersion of a
surface in a 3--manifold and $f_*(\pi_1(\Sigma))$ is a separable
subgroup of $\pi_1(M)$ then the immersion factors (up to homotopy)
through an embedding in a finite cover of $M$. This technique has
applications to the still open ``virtual Haken conjecture'' and the
``positive virtual first Betti number conjecture''.

\proclaim{The virtual Haken conjecture} If $M$ is a compact,
irreducible 3--manifold with infinite fundamental group then $M$ is virtually
Haken, that is it has a finite cover which contains an embedded, 2--sided,
incompressible surface.
\endproc

\proclaim{The positive virtual first Betti number conjecture}If $M$ is a
compact,
irreducible 3--manifold with infinite fundamental group then it has
a finite cover with positive first Betti number.
\endproc

Unfortunately it is difficult in
general to show that a given subgroup is separable, and it is known that not
every
subgroup of a 3--manifold group need be separable; the first example was
given by
Burns, Karrass and Solitar, [\BKS]. On the other hand Shalen has shown that
if an
aspherical 3--manifold admits a
$\pi_1$--injective immersion of a surface which factors through infinitely
many
finite covers then the 3--manifold is virtually Haken [\N]. In group
theoretic
terms Shalen's condition says that the surface subgroup is
contained in infinitely many finite index subgroups of the fundamental group
of the 3--manifold, and this is clearly a weaker requirement than 
separability.

The engulfing property is apparently weaker still. It was introduced by Long
in
[\Long] to study hyperbolic 3--manifolds, and he was able to show that in
some
circumstances it implies separability.  He remarks that ``One of the
difficulties with the LERF (separability) property is that there often
appears
to be nowhere to start, that is, it is conceivable that a finitely generated
proper subgroup could be contained in no proper subgroups of finite index at
all.'' In this note we show that this can happen  for  finitely
generated subgroups of the fundamental group of a Haken (though not
hyperbolic)
3--manifold. We give  two examples, both already known not to be subgroup
separable. One is derived from the  recent work of Rubinstein and Wang,
[\RW],
and we consider it in Theorem 1. The other was the first known example of a
3--manifold group which failed to be subgroup separable and was introduced in
[\BKS] and further studied in [\LN] and [\NW]. Our proof that it fails to
satisfy the engulfing property is more elementary than the original proof
that it
fails subgroup separability, and we hope that it sheds
some light on this fact. Both of the examples are graph manifolds so
they leave open the question of whether or not hyperbolic 3--manifold groups 
are
subgroup separable or satisfy the engulfing property.  In this connection we
note that if every surface subgroup of  any closed hyperbolic 3--manifold
does
satisfy the engulfing property then any such subgroup must be contained in
infinitely many finite index subgroups, and Shalen's theorem would give a
solution to the ``virtual Haken conjecture'' for closed hyperbolic
3--manifolds
containing surface subgroups.

\section{The example of Rubinstein and Wang}

We will use the following lemma:
\proclaim{Lemma 1} Let $H$ be a separable subgroup of a group $G$. Then
the index $[G:H]$ is finite if and only if there is a finite subset
$F\subset G$ such that $G=HFH$.
\endproc

\prf If $[G:H]$ is finite then $G=FH$ for some finite subset
$F\subseteq G$, so $G=HFH$ as required.

Now suppose that $G=HFH$ for some finite subset $F\subset G$. For each
element $g\in F\setminus(H\cap F)$, we can find a finite index subgroup
$H_g\in G$ with $H<H_g$ but $g\not \in H_g$. Now let $K=\mathop\cap\limits_g
H_g$. Since $F$ is finite, $K$ has finite index in $G$, and since
$H<K$, $K$ contains every double coset $HgH$ which it intersects
non-trivially. It follows that $K$ only intersects a double coset $HgH$
non-trivially if $g\in H$, and so $K=H$.
\endprf

Given a subgroup $H<G$ let $\overline H$ denote the intersection of the
finite index subgroups of $G$ which contain $H$. ($\overline H$ is the 
closure
of $H$ in the profinite topology on $G$). It is obvious that $H$ is
separable
if and only if $H=\overline H$, and it is engulfed if and only if $G\neq
\overline H$. If
$G$ is a finite union of double cosets of a subgroup $H$ then it is also a
finite
union of double cosets of $\overline H$ and this is clearly a separable
subgroup of $G$ so by Lemma 1 it must have
finite index. Now if $H$ has infinite index in $G$ and $\overline H$ has
finite index in $G$ they cannot be equal, and $H$ is not separable. Hence we
may interpret a finite double coset decomposition
$G=HFH$ as an obstruction to separability for an infinite index subgroup
$H<G$.

In [\RW] Rubinstein and Wang constructed a graph manifold $M$ and a
$\pi_1$--injective immersion $\phi\co \Sigma\looparrowright M$ of a surface 
$\Sigma$ which does not factor through an embedding into any
finite cover of $M$. It follows from [\Sc] that the surface group
$H=\phi_*(\pi_1(\Sigma))$ is not separable in the 3--manifold group
$G=\pi_1(M)$.   In fact as we shall see  $G$ has a finite double
coset decomposition $G=HFH$:

\proclaim{Lemma 2} Let $\phi\co \Sigma\looparrowright M$ be a
$\pi_1$--injective immersion of a surface $\Sigma$ in a 3--manifold $M$,
and let $M_H$ be the cover of $M$ defined by the inclusion
$\phi_*(\pi_1(\Sigma))\hookrightarrow\pi_1(M)$. Let
$\tilde{\phi}\co \R^2\looparrowright\tilde{M}$ be some lift of $\phi$ to the
universal covers, and $\tilde{\Sigma}$ denote the image
of
$\tilde{\phi}$. Then the number of $H$ orbits for the action on
$G\tilde{\Sigma}=\{g\tilde{\Sigma}\mid g\in G\}$ is precisely the
number of distinct double cosets
$HgH$.
\endproc

\prf By construction $\tilde{\Sigma}$ is $H$--invariant, so for each
double coset $HgH$ we have
$HgH\tilde{\Sigma}=Hg\tilde{\Sigma}$. It follows that if $F=\{g_i\mid i\in
I\}$ is a complete family of representatives for the
distinct double cosets
$Hg_iH$ in $G$ then the $G$--orbit $G\tilde{\Sigma}$ breaks into $|F|$
$H$--orbits as required.
\endprf

Now in the example in [\RW] we are told in Corollary 2.5 that the image of
each orbit $Hg\tilde(\Sigma)$  intersects the image of
$H\tilde{\Sigma}$ which by construction of $H$ is compact. Hence there are
only finitely many such images, and therefore only
finitely many $H$--orbits for the action of $H$ on the set
$G\tilde{\Sigma}$. Hence $G=HFH$ for some finite subset $F\subset G$.

\proclaim{Corollary} The profinite closure of
$H$ must have finite index in
$G$, ie there are only finitely many finite index subgroups of $G$
containing $H$, or, in
topological terms, there are only finitely many finite covers of the
3--manifold $M$ to which the
surface $\Sigma$ lifts by degree~1.
\endproc

Now as in the proof of Lemma 1, let $K$ denote the intersection of the
finite index subgroups of $G$ containing $H$, and let  $M_K$
denote the finite cover of
$M$ corresponding to the finite index subgroup $K<G$. Then the immersion of
$\Sigma$ in $M$
lifts to an immersion
$\bar\phi\co \Sigma\looparrowright M_K$ which does not lift to any
finite cover of
$M_K$. Hence:

\proclaim{Theorem 1} There is a compact 3--manifold $M_K$ and a
$\pi_1$--injective
immersion
$\bar\phi\co \Sigma\looparrowright M_K$ which does not factor through any
proper finite
cover of
$M_K$.
\endproc

\section{The example of Burns, Karrass and Solitar}

In [\BKS], Burns Karrass and Solitar gave an
example of a 3--manifold group with a finitely generated subgroup
which is not separable.
Their example is a free by $\Z$ group with presentation
$\langle  \alpha, \beta, y \mid \alpha^y=\alpha\beta, \beta^y=\beta\rangle$.
It is easy to show that their example is isomorphic
to the group $G$ with presentation
$\langle a,b,t \mid [a,b], a^t=b \rangle$,
and it is in this form that we shall work with $G$.
Note that here and below we use the notation
$x^y = y^{-1}xy$ and $[x,y]= x^{-1}y^{1}xy$.

In this section we show that $G$ has a proper subgroup
$K\subset G$ such that $K$
is not engulfed.
 In particular, this yields an easier
proof that $G$ has non-separable subgroups.

\proclaim{Lemma 3}
Let $J=\langle abb, t \rangle$.
Let $H$ be a finite index subgroup of $G$ containing $J$.
Then $G = H\langle a\rangle$.
\endproc

\prf
We express the argument in terms of covering spaces.
Let $X$ denote the standard based 2--complex for the presentation of
$G$.
Let $T$ denote the torus subcomplex $\langle a,b \mid [a,b] \rangle$
of $X$.
The complex $X$ is formed from $T$ by the addition of a cylinder
$C$ whose top and bottom boundary components are attached to the
loops $a$ and $b$ respectively, and $C$ is subdivided by
a single edge labeled~$t$ which is oriented from the
$a$ loop to the $b$ loop.

Let $\hat X$ denote the finite based cover of $X$ corresponding to the
subgroup $H$.
Let $\hat T$ denote the cover of $T$ at the basepoint of $\hat X$.
Let $\hat a$ and $\hat b$ denote the covers of the loops $a$ and $b$
at the basepoint.

Since  $t$ lifts to a closed path
in $\hat X$, we see that $C$ has a finite cover $\hat C$ which
lifts at the basepoint
to a cylinder attached at  its ends to $\hat a$ and
$\hat b$.
Now $\hat C$
gives a one-to-one correspondence between $0$--cells on
$\hat a$ and $0$--cells on $\hat b$.
In particular, each $t$~edge of $\hat C$ is directed
from some $0$--cell in $\hat a$ to some $0$--cell in $\hat b$ and therefore
 $\deg{\hat a} = \deg{\hat b}$.

Because $abb \in J\subset H$ and hence $abb \in \pi_1(\hat T)$,
we see that $b$  generates the covering group of the regular cover
$\hat{T}\longrightarrow T$,
and therefore $\deg{\hat b} = \deg{\hat T}$.
Thus we have $\deg{\hat T} = \deg{\hat b} = \deg{\hat a}$,
and because $\deg{\hat T}$ is finite,
we see that every $0$--cell of $\hat T$ lies in both $\hat a$
and  $\hat b$.

As above, each $0$--cell of $\hat a$ has an outgoing $t$~edge
in $\hat C$ and each $0$--cell of $\hat b$ has an incoming $t$~edge
in $\hat C$, and so we see that each $0$--cell of $\hat T \cup \hat C$ has
an incoming
and outgoing $t$~edge.
Since $0$--cells of $\hat T \cup \hat C$ obviously have incoming and outgoing
$a$ and $b$ edges in $\hat T$,
we see that $\hat X = \hat T \cup \hat C$ and in particular,
every $0$--cell of $\hat X$ is contained in $\hat T$ and therefore in $\hat
a$.
Thus  $\langle a \rangle$ contains a
set of right coset representatives for $H$ in $G$,
and consequently $G = H\langle a\rangle$.
\endprf

\proclaim{Lemma 4}
Let $K = \langle J \cup a^g\rangle$ for some $g\in G$.
Then $K$ is not engulfed.
\endproc

\prf
Let $H$ be a subgroup of finite index containing $K$.
Since $J\subset H$ we may apply Lemma~3 to conclude that
$G = H\langle a\rangle$ and so it is sufficient to show that
$a \in H$.
Observe that $g^{-1}=ha^n$ for some $h\in H$ and $n\in \Z$.
But $a^g = (ha^n)aa^{-n}h^{-1} = hah^{-1}$,
and obviously $hah^{-1} \in H$ implies that $a \in H$.\break
\endprf

\proclaim{Theorem~2}
Let $K$ be the subgroup $\langle abb, t, btat^{-1}b^{-1} \rangle$.
Then the engulfing property fails for $K$, that is, $K\neq G$ and
the only subgroup of finite index containing $K$ is $G$.
\endproc

\prf
Lemma~4 with $g=t^{-1}b^{-1}$ shows that
$K$ is not engulfed.
To see that $K \neq G$ we observe  that the normal
form theorem for an HNN extension shows that there is no non-trivial
cancellation
between the generators of $K$ so it is a rank~3 free group, but $G$ is not
free.
\endprf

\rk{Remark}
\rm It is not difficult to see that there are many  finitely generated
subgroups   $J$ for which some
version of Lemma~3 is true. In addition, one has some freedom to
vary the choice of $g$ in theorem~2.
Consequently  subgroups of $G$  which are not engulfed are numerous.

\references     
\Addresses\recd
\bye